\documentclass[a4paper]{amsart}
\usepackage[british]{babel}
\usepackage{amsmath}
\usepackage{amsthm}
\usepackage{amscd}
\usepackage{amsfonts}
\usepackage{amssymb}
\usepackage{hyperref}
\usepackage{ifthen}

\newcommand{\Q}{\mathbb{Q}}
\newcommand{\Z}{\mathbb{Z}}
\newcommand{\X}{\mathcal{X}}
\newcommand{\Xs}{{\mathcal{X}_s}}
\newcommand{\Ps}{P_s}
\newcommand{\Xb}{\bar{X}}
\newcommand{\Xsb}{\bar{\mathcal{X}}_s}
\newcommand{\G}{\mathcal{G}}
\newcommand{\Gs}{\mathcal{G}_s}
\newcommand{\Gm}{\mathbb{G}_\textrm{m}}

\renewcommand{\H}{\mathrm{H}}
\newcommand{\R}{\mathrm{R}}
\renewcommand{\O}{\mathcal{O}}

\newcommand{\Knr}{{K^\textrm{nr}}}
\newcommand{\Ksep}{K^s}

\newcommand{\A}{\mathbb{A}}
\newcommand{\T}{\mathcal{T}}

\DeclareMathOperator{\Pic}{Pic}
\DeclareMathOperator{\Spec}{Spec}
\DeclareMathOperator{\Gal}{Gal}
\DeclareMathOperator{\Hom}{Hom}
\newcommand{\kb}{{\bar{k}}}
\newcommand{\et}{\mathrm{\acute{e}t}}
\newcommand{\sm}{\mathrm{sm}}
\newcommand{\WR}{\mathbb{R}}

\newtheorem{theorem}{Theorem}[section]
\newtheorem{lemma}[theorem]{Lemma}
\newtheorem{proposition}[theorem]{Proposition}
\newtheorem{corollary}[theorem]{Corollary}
\theoremstyle{remark}
\newtheorem*{remark}{Remark}

\theoremstyle{definition}

\author{Martin Bright}
\title{Torsors under tori and N\'eron models}
\email{M.Bright@warwick.ac.uk}
\address{Martin Bright \\ Mathematics Institute \\ Zeeman Building \\ 
University of Warwick \\ Coventry CV4 7AL}
\thanks{This research was funded by the Heilbronn Institute for
Mathematical Research.}

\begin{document}


\begin{abstract}
Let $R$ be a Henselian discrete valuation ring with field of fractions
$K$.  If $X$ is a smooth variety over $K$ and $G$ a torus over $K$,
then we consider $X$-torsors under $G$.  If $\X/R$ is a model of $X$
then, using a result of Brahm, we show that $X$-torsors under $G$
extend to $\X$-torsors under a N\'eron model of $G$ if $G$ is split
by a tamely ramified extension of $K$.  It follows that
the evaluation map associated to such a torsor factors through
reduction to the special fibre.  In this way we can use the geometry
of the special fibre to study the arithmetic of $X$.
\end{abstract}

\maketitle

\section{Introduction}

Let $R$ be a Henselian discrete valuation ring, with field of
fractions $K$ and perfect residue field $k$. Let $s \colon \Spec k \to
\Spec R$ denote the inclusion of the special point.

If $X$ is a smooth variety\footnote{By a variety over $K$ we mean a
  separated, integral scheme of finite type over $K$.} over $K$, and
$Y$ an $X$-torsor under an algebraic $K$-group $G$, there is an
``evaluation'' map $X(K) \to \H^1(K,G)$ which associates to each point
$P$ of $X$ the isomorphism class of the fibre $Y_P$.  The main result
of this article is to show that, when $G$ is a torus split by a tamely
ramified extension of $K$, this map depends only on the image of $P$
in the special fibre of a model of $X$.

\begin{theorem}\label{thm:main}
Let $X$ be a smooth variety over $K$, and let $\X/R$ be a model of
$X$; we write $X(K)_\sm$ for the set of $K$-points of $X$ which extend
to smooth $R$-points of $\X$.  Let $G$ be a torus over $K$ with a
N\'eron model $\G$ over $R$, and suppose that $G$ is split by a tamely
ramified extension of $K$.  Let $Y \to X$ be an $X$-torsor under $G$.
Then the natural maps
\[
\H^1(K,G) \leftarrow \H^1(R,\G) \to \H^1(k,\Gs)
\]
are isomorphisms, and the evaluation map $X(K)_\sm \to \H^1(K,G)$
coming from $Y$ factors through the $k$-points of the special fibre
$\Xs$, as follows:
\[
\begin{CD}
X(K)_\sm @>>> \Xs(k) \\
@V{Y}VV     @VV{f}V \\
\H^1(K,G) @>{\cong}>> \H^1(k,\Gs) \text{.}
\end{CD}
\]
Here $\Gs$ is the special fibre of $\G$, and the map $f$ comes from an
element of $\H^1(\Xs, \Gs)$.
\end{theorem}

A main step in proving Theorem~\ref{thm:main} is the following
technical result, in which the principal ingredient is due to Brahm.

\begin{theorem}\label{thm:tech}
Let $\X$ be a smooth scheme over $R$.  Let $G$ be a torus over $K$ with a
N\'eron model $\G$ over $R$, and suppose that $G$ is split by a tamely
ramified extension of $K$.  Then the natural map $\H^1(\X,\G) \to
\H^1(X,G)$ is an isomorphism.
\end{theorem}

If the residue field $k$ is in fact finite, then we can easily deduce
the following refinement of Theorem~\ref{thm:main}, which is more
useful in applications.
                                 
\begin{corollary} \label{cor:finite}
Suppose in addition that $k$ is finite.  Let $\Phi(\G)$ denote the
$k$-group of connected components of $\G$.  Then we can replace $\Gs$
with $\Phi(\G)$ in the conclusion of Theorem~\ref{thm:main}: the
evaluation map $X(K)_\sm \to \H^1(K,G)$ coming from $Y$ factors as follows:
\[
\begin{CD}
X(K)_\sm @>>> \Xs(k) \\
@V{Y}VV     @VV{g}V \\
\H^1(K,G) @>{\cong}>> \H^1(k,\Phi(\G)) \text{.}
\end{CD}
\]
Here the map $g$ comes from an element of $\H^1(\Xs, \Phi(\G))$.
\end{corollary}

This form of the theorem is particularly useful because $\H^1(\Xs,
\Phi(\G))$ is often easy to calculate, if the geometry of $\Xs$ is
known.

\begin{remark}
If the model $\X$ is proper over $R$ and regular, then in fact every
$K$-point of $X$ extends to a smooth $R$-point of $\X$:
see~\cite[Chapter~3, Proposition~2]{BLR:NM}.
\end{remark}

Our proof of Theorem~\ref{thm:tech}, and hence Theorem~\ref{thm:main},
will in fact use only the following properties of the torus $G$:
\begin{enumerate}
\item \label{hyp1} $G$ is a smooth, commutative algebraic group over
  $K$ admitting a N\'eron model over $R$;
\item \label{hyp2} if $\eta \colon \Spec K \to \Spec R$ denotes the
  inclusion of the generic point, then $\R^1 \eta_* G = 0$, considered
  as a sheaf on the smooth site over $R$.
\end{enumerate}
That~(\ref{hyp2}) holds for any torus split by a tamely ramified
extension was proved by Brahm~\cite[Chapter~4]{Brahm:SMIUM-2004}, and
relies on our assumption that the residue field is perfect.  Brahm
also studied the structure of $\R^1 \eta_* G$ when this assumption is
removed.

\subsection{Motivation}

Recall that a $K$-torus is an affine group scheme over $K$ which
becomes isomorphic to a product of some number of copies of $\Gm$
after a finite base extension; we say that such an extension
\emph{splits} the torus.  As described in~\cite[Chapter~10]{BLR:NM},
tori have N\'eron models which are locally of finite type, but may not
be of finite type.  The group of components of the N\'eron model is a
finitely generated Abelian group with an action of $\Gal(\kb/k)$; it
has been described by Xarles~\cite{Xarles:JRAM-1993}.

For example, let $K=\Q_p$ with $p$ odd, let $L=\Q_p(\sqrt{p})$, a
tamely ramified quadratic extension of $K$, and let $T = \WR^1_{L/K}
\Gm$ be the norm torus associated to this extension.  The component
group of the N\'eron model of $T$ is of order $2$.  It is well known
that $\H^1(K,T) = K^\times / N L^\times$, and that there an
isomorphism of this group with $k^\times / (k^\times)^2 = \H^1(k,
\Z/2\Z)$; see~\cite[Chapter~V, \S 3]{Serre:CL}.  For units in $K$, the
isomorphism is simply reduction modulo $p$: $\O_K^\times / N_{L/K}
\O_L^\times \cong k^\times / (k^\times)^2$.

Let $X$ be a smooth variety over $K=\Q_p$, and $Y$ an $X$-torsor under
$T$.  On a sufficiently small open subsets $U \subset X$, $Y$ looks
like
\[
\{ x^2 - p y^2 = f \} \subset \A^2_U
\]
for some $f \in K[U]^\times$.  For any $P \in U(K)$, the isomorphism
class of the fibre $Y_P$ is thus given by the class of $f(P)$ in
$K^\times/ NL^\times$.  If $f$ takes values in $\Z_p^\times$ on a
neighbourhood of $P$, then we see that in fact the isomorphism class
of $Y_P$ depends only on whether the reduction modulo $p$ of $f(P)$ is
a square, and hence depends only on the reduction of $P$ modulo $p$.

Given a model $\X/\Z_p$ of $X$, it is not too hard to show that,
whenever $P$ extends to a smooth $\Z_p$-point of $\X$, we can arrange
for $f$ to take values in $\Z_p^\times$ on a neighbourhood of $P$.  We
have therefore shown that there is a commutative diagram
\[
\begin{CD}
X(K)_\sm @>>> \Xs(k) \\
@V{Y}VV     @VVV \\
\H^1(K,T) = K^\times / NL^\times @>{\cong}>> 
k^\times/(k^\times)^2 = \H^1(k,\Z/2\Z) \text{.}
\end{CD}
\]

The aim of this article is to generalise this example.  In doing so,
there are two principal questions:
\begin{itemize}
\item What other algebraic groups may replace the torus $T$?
\item Given a suitable algebraic group, what should replace $\Z/2\Z$
  in this example?
\end{itemize}

Theorem~\ref{thm:main} states that $T$ may be replaced by any
$K$-torus which is split by a tamely ramified extension, and that
$\Z/2\Z$ should be replaced by the component group of the
corresponding N\'eron model.

\begin{remark}
It is not reasonable to expect Theorem~\ref{thm:main} to hold for all
torsors under tori.  Indeed, consider the norm torus $T =
\WR^1_{L/\Q_2} \Gm$ for a quadratic, ramified extension $L/\Q_2$.  The
question of whether $a \in \Z_2^\times$ is a norm from $L$ does not
depend only on $a$ modulo $2$: rather, it depends on $a$ modulo $4$ or
$8$.  We cannot therefore expect the evaluation map coming from a
torsor under $T$ to factor through reduction modulo $2$.
\end{remark}

\subsection{Background}

The use of torsors to study the rational points of varieties has been
an important part of the theory of Diophantine equations since the
middle of the 20th century; an excellent reference for torsors and
their applications to the study of rational points
is~\cite{Skorobogatov:TRP}.  In the study of rational varieties,
Colliot-Th\'el\`ene and Sansuc showed that the Brauer--Manin
obstruction is naturally studied using torsors under
tori~\cite{CTS:DMJ-1987}; indeed, they conjectured that the
obstructions coming from torsors under tori are the only obstructions,
both to the existence of rational points and to weak approximation, on
rational varieties.  It is therefore important to understand the
arithmetic and geometry of torsors under tori for two reasons: to
approach this conjecture; and to develop ways of calculating the
obstructions.  The object of this article is to describe one way in
which the arithmetic of torsors on a variety $X$ is linked to the
geometry of models of $X$, thus giving a new way to study them.

In this article, we concern ourselves entirely with the local picture:
we fix a local field $K$, a smooth variety $X$ over $K$, and a torsor
$Y \to X$ under some $K$-group $G$.  We show that, under certain
hypotheses on $G$, the evaluation map $X(K) \to \H^1(K,G)$ factors
through reduction to the special fibre: that is, for a suitable model
$\X$ of $X$, the evaluation map comes from a map $\X(k) \to \H^1(k,G)$,
where $k$ is the residue field of $K$ and the map comes from a torsor
under a $k$-group.  In this way, the arithmetic of the torsor $Y$ is
related to the geometry of the special fibre of $\X$; we give examples
of how this can usefully be used to deduce facts about $Y$.

Although our principal application is to torsors under tori, many of
the results collected will apply equally well to torsors under more
general algebraic groups, and we will state them in this generality
where possible.

The method used to obtain these results is to look at a N\'eron model
of the group $G$, and try to extend $X$-torsors under $G$ to torsors
over $\X$ under a N\'eron model of $G$, if one exists.  In
Section~\ref{sec:constant}, we collect some results about constant
torsors: how do $K$-torsors under a group $G$ extend to torsors under
a N\'eron model?  In Section~\ref{sec:varieties}, we ask the same
question for $X$-torsors, and deduce the main result
(Theorem~\ref{thm:main}).

\subsection{Notation}

Throughout the rest of this article, $R$ denotes a Henselian discrete
valuation ring, with field of fractions $K$ and perfect (but not
necessarily finite) residue field $k$.  Let $\Ksep$ denote a separable
closure of $K$, and $\kb$ the corresponding algebraic closure of $k$.
Let $\Knr$ denote the maximal unramified extension of $K$.  Let $\eta
\colon \Spec K \to \Spec R$ be the inclusion of the generic point, and
let $s \colon \Spec k \to \Spec R$ be the inclusion of the special
point.

If $S$ is a scheme, then the site $S_\et$ is the small \'etale site,
the category of schemes \'etale and of finite type over $S$, endowed
with the \'etale topology, for which the coverings are surjective
families of \'etale morphisms of finite type.  The site $S_\sm$ is the
smooth site, the category of schemes smooth over $S$ with the topology
in which coverings are surjective families of smooth morphisms.

Suppose that $X$ is a smooth variety over $K$.  By \emph{model} of $X$
over $R$ we mean a scheme $\X$, separated, flat and locally of finite
type over $R$, with generic fibre $\X_\eta$ isomorphic to $X$.
Suppose that $\X/R$ is a model of $X$.  If $\X$ is not proper, then
$K$-points of $X$ do not necessarily extend to $R$-points of $\X$;
even if $\X$ is proper, then $K$-points of $X$ cannot be expected to
extend to \emph{smooth} $R$-points of $\X$.  In either situation,
define $X(K)_\sm$ to be the set of $K$-points of $X$ which do extend
to smooth points of $\X$.

If $G$ is a smooth group scheme over $K$, then a \emph{N\'eron model}
of $G$ is a smooth and separated group scheme, locally of finite type,
over $R$ which represents the sheaf $\eta_* G$ on the smooth site
$R_\sm$.  N\'eron's original construction applied to Abelian
varieties, but Raynaud showed that N\'eron models also exist for some
affine algebraic groups, and in particular for
tori~\cite[Chapter~10]{BLR:NM}.  Given a group scheme $\G$ over $R$,
let $\Phi(\G)$ denote the group of connected components of $\G$,
considered either as a group scheme over $k$ or as a
$\Gal(\kb/k)$-module.

Cohomology groups, unless otherwise indicated, are \'etale cohomology
groups.  Since we will be dealing only with smooth group schemes, the
\'etale topology is fine enough to classify torsors --
see~\cite[Chapter~III, Remark~4.8]{Milne:EC}.

If $S$ is any scheme and $\G$ a smooth, commutative group scheme over
$S$, then $\G$ gives rise to a sheaf on the \'etale site $X_\et$,
where $X$ is any scheme locally of finite type over $S$.  (Indeed,
$\G$ defines a sheaf on the big \'etale site over $S$.)  If $Y \to X$
is any morphism of $S$-schemes, then there is a natural morphism
$\G(X) \to \G(Y)$, functorial in $\G$; by the universal
$\delta$-functor property of cohomology groups, this gives rise to a
natural morphism $\H^i(X,\G) \to \H^i(Y,\G)$ for each $i \ge 0$, again
functorial in $\G$.  In particular, let $\alpha \in \H^1(X,\G)$ be
fixed, and let the morphism $Y \to X$ vary; we obtain a map $X(Y) \to
\H^1(Y,\G)$, which we will denote again by $\alpha$, and call it the
\emph{evaluation map} associated to $\alpha$.

Let $P \to X$ be an $X$-torsor under $\G$; then we can associate to
$P$ a class $\alpha$ in $\H^1(X,\G)$.  The pull-back $P \times_X Y$ is
a $Y$-torsor under $\G$, and the corresponding class in $\H^1(Y,\G)$
is the image of $\alpha$ under the natural morphism just defined.
See~\cite[Chapitre~III, \S 2.4.6.2]{Giraud:CNA}.  In this way the
evaluation map $X(Y) \to \H^1(Y,\G)$ is seen to be the map which
associates to each point $y \in X(Y)$ the isomorphism class of the
fibre $P_y$.

\section{Extending torsors}

Let $G$ be a smooth algebraic group over $K$, with a N\'eron model
$\G$.  Our aim is to understand how torsors under $G$ over some
$K$-scheme extend to torsors under $\G$.  In this section we collect
some known results and show how they answer this question.

\subsection{Constant torsors} \label{sec:constant}
If we are to have any hope of understanding torsors on varieties, we
should first look at the question for constant torsors: how do
$K$-torsors under $G$ extend to $R$-torsors under $\G$?  This is
essentially a question in Galois cohomology, but we present it using
\'etale cohomology for consistency with the next section.

In fact we will take a purely cohomological view of this problem, and
ask about extending classes in $\H^1(K,G)$ to classes in
$\H^1(R,\G)$.  Whether the resulting cohomology classes are actually
represented by torsors is unimportant for our purposes.

The facts presented in this section are no more difficult to state for
non-commutative groups $G$, and so we remind the reader of the
definitions of non-commutative $\H^1$ in Galois
cohomology~\cite[Chapitre~1, \S 5]{Serre:CG} and in \'etale
cohomology~\cite[Chapitre~III]{Giraud:CNA} or~\cite[III.4]{Milne:EC}.
We will also need the definition of the non-commutative higher direct
image $\R^1$: see~\cite[Chapitre~V, \S 2]{Giraud:CNA}.

We begin by recalling that, for any connected linear algebraic group
$G$ over $K$, we have $\H^1(\Knr,G)=0$ by Th\'eor\`eme~$1'$
of~\cite[Chapitre~III, \S 2.3]{Serre:CG}.  (This result depends on the
fact that $\Knr$ is a field of dimension $\leq 1$, which in turn
relies on our hypothesis that the residue field $k$ be perfect.)

\begin{lemma}\label{lem:R1et}
Let $G$ be a smooth algebraic group over $K$.  Suppose that
$\H^1(\Knr,G)=0$.  Then $\R^1\eta_* G = 0$ as a sheaf on $R_\et$.
\end{lemma}
\begin{proof}
The \'etale sheaf $\R^1 \eta_* G$ has two stalks: that at the generic
point is $\H^1(\Ksep,G)$, which is trivial since $\Ksep$ is separably
closed; that at the special point is $\H^1(\Knr,G)$, which is trivial
by assumption.  Therefore $\R^1 \eta_* G = 0$ on $R_\et$.
\end{proof}

\begin{lemma}\label{lem:const}
Let $G$ be a smooth algebraic group over $K$; suppose that $G$ admits
a N\'eron model $\G$, and suppose that $\H^1(\Knr,G)=0$.  Then there
is an isomorphism (of pointed sets, or of groups if $G$ is
commutative)
\begin{equation}
\H^1(K,G) \cong \H^1(R, \G) \text{.}
\end{equation}
\end{lemma}
\begin{proof}
If $G$ is commutative, this follows from the Leray spectral sequence
for $\eta$.  More generally, there is an exact sequence
\cite[Chapitre~V, Proposition~3.1.3]{Giraud:CNA} as follows
\[
1 \to \H^1(R, \eta_* G) \xrightarrow{\alpha} \H^1(K, G) \to \H^0(R,
\R^1 \eta_* G) \text{.}
\]
Although this is only an exact sequence of pointed sets, twisting
shows that $\alpha$ is injective.  Since $\G$ represents the sheaf
$\eta_* G$ on $R_\sm$ and hence on $R_\et$, Lemma~\ref{lem:R1et} shows
that $\alpha$ gives the desired isomorphism.
\end{proof}

\begin{remark}
This is really the same as using the inflation-restriction sequence in
Galois cohomology: from the exact sequence
\begin{multline*}
0 \to \H^1(\Gal(\Knr/K), G(\Knr)) \xrightarrow{\text{inf}}
\H^1(\Gal(\Ksep/K), G(\Ksep)) \\
\xrightarrow{\text{res}} \H^1(\Gal(\Ksep/\Knr), G(\Knr))
\end{multline*}
we deduce an isomorphism 
\[
\H^1(\Gal(\Ksep/K),G(\Ksep)) \cong \H^1(\Gal(\Knr/K), G(\Knr)).
\]
The first of these groups is the same as $\H^1(K,G)$, using the usual
correspondence between sheaves on $K_\et$ and $\Gal(\Ksep/K)$-sets.
The second is isomorphic to $\H^1(R,\G)$ since $\G(R^\textrm{nr}) =
G(\Knr)$, where $R^\textrm{nr}$ is the integral closure of $R$ in
$\Knr$, and \'etale cohomology of a Henselian local ring may be
computed using finite \'etale covers.
\end{remark}

\begin{lemma}\label{lem:hensel}
Let $\G$ be a smooth group scheme over $R$, not necessarily of finite
type.  Then there is an isomorphism
\begin{equation}
\H^1(R,\G) \cong \H^1(k, \Gs) \text{.}
\end{equation}
\end{lemma}
\begin{proof}
That the natural map $\H^1(R,\G) \to \H^1(k, \Gs)$ is an isomorphism
is well known, and essentially due to Hensel's Lemma, but there are
some issues around representability to be dealt with before Hensel's
Lemma can be applied.  A very general proof of this fact, which does
not require that $\G$ be of finite type over $R$, can be found
at~\cite[Th\'eor\`eme~11.7 and Remarque~11.8(3)]{Grothendieck:GB3}.
\end{proof}

When the residue field $k$ is finite, we can go one step further.
\begin{lemma}\label{lem:finite}
Suppose further that $k$ is finite, and that $\G$ is either commutative
or of finite type over $k$.  Then there is an isomorphism
\[
\H^1(k, \Gs) \to \H^1(k, \Phi(\G)) \text{.}
\]
\end{lemma}

\begin{proof}
Consider the exact sequence of group schemes over $k$:
\[
0 \to (\Gs)_0 \to \Gs \to \Phi(\G) \to 0
\]
where $(\Gs)_0$ is the connected component of $\Gs$ containing the
identity.  If $G$ is commutative, the theorem of
Lang~\cite[Chapter~VI, Proposition~5]{Serre:AGCF} shows that $\H^i(k,
(\Gs)_0)=0$ for $i>0$, and so $\H^1(k, \Gs) \to \H^1(k, \Phi(\G))$ is
an isomorphism.  If $G$ is not commutative, the same theorem of Lang
still gives $\H^1(k, (\Gs)_0)=0$, showing that our map is injective.
When $\G$ is of finite type, then $\Gs$ is an algebraic group over $k$
and Corollary~3 of~\cite[Chapitre~III, \S 2.4]{Serre:CG} shows that
the map is also surjective.
\end{proof}

\subsection{Torsors over varieties}\label{sec:varieties}

Let $X$ be a variety over $K$.  In this theorem, we describe when a
torsor $Y \to X$ under a group $G$ extends to a (sheaf) torsor under
the N\'eron model of $G$ over a smooth model of $X$.  When this does
happen, we can deduce that the evaluation map associated to $Y$
factors through the $k$-points of the model, resulting in
Theorem~\ref{thm:main}.

\begin{lemma}\label{lem:extension}
Let $\X$ be a smooth scheme over $R$.  Let $G$ be a smooth algebraic
group over $K$ with a N\'eron model $\G$ over $R$, and suppose that
$\R^1 \eta_* G = 0$ as a sheaf on the smooth site $R_\sm$.  Then the
natural map $\H^1_\et(\X, \G) \to \H^1_\et(\X_\eta, G)$ is an
isomorphism.
\end{lemma}
\begin{proof}
Firstly note that, since $\G$ is smooth over $R$, both of
$\H^1_\et(\X, \G)$ and $\H^1_\et(\X_\eta, G)$ are isomorphic to the
corresponding cohomology groups or sets on the respective smooth
sites.  (For a proof, see~\cite[XXIV, Proposition~8.1]{SGA3III}
or~\cite[Remarque~11.8(3)]{Grothendieck:GB3}.)

Now apply~\cite[Chapitre~V, Proposition~3.1.3]{Giraud:CNA} to the
inclusion of the generic fibre $j\colon (\X_\eta)_\sm \to \X_\sm$, to
obtain an exact sequence of pointed sets
\[
1 \to \H^1_\sm(\X, j_* G) \xrightarrow{\alpha} \H^1_\sm(\X_\eta,G) \to
\H^0_\sm(\X, \R^1 j_* G) \text{.}
\]
Since $\X$ is smooth over $R$, the sheaf $j_* G$ is represented by
$\G$, and the sheaf $\R^1 j_* G$ on $\X_\sm$ is simply the restriction
to $\X_\sm$ of the sheaf $\R^1 \eta_* G$ on $R_\sm$.  The hypothesis
implies that $\R^1 j_* G = 0$ on $\X_\sm$, and so $\alpha$ is
surjective.  From the exact sequence above, a twisting argument shows
that $\alpha$ is also injective.
\end{proof}

In particular, we can deduce Theorem~\ref{thm:tech}.

\begin{proof}[Proof of Theorem~\ref{thm:tech}]
Brahm~\cite[Chapter~4]{Brahm:SMIUM-2004} has shown that a torus split
by a tamely ramified extension satisfies the hypotheses of
Lemma~\ref{lem:extension}.
\end{proof}

Note that, in contrast to the situation for constant torsors, we can
in general say nothing interesting about the natural map from
$\H^1(\X, \G)$ to $\H^1(\Xs,\Gs)$.

Let us now derive Theorem~\ref{thm:main} by putting
Theorem~\ref{thm:tech} together with the facts collected in
Section~\ref{sec:constant}.

\begin{lemma}\label{lem:eval}
Let $\X$ be a smooth scheme over $R$, $\G$ a group scheme over $R$,
and $P \in \X(R)$ a point of $\X$.  Then the following diagram of
evaluation maps commutes:
\[
\begin{CD}
\H^1(\Xs, \Phi(\G)) @<<< \H^1(\Xs, \Gs) @<<< \H^1(\X,\G) 
                                        @>>> \H^1(\X_\eta, \G_\eta) \\
@VV{\Ps}V                @VV{\Ps}V           @VV{P}V          
                                             @VV{P_\eta}V        \\
\H^1(k, \Phi(\G))   @<<< \H^1(k, \Gs)   @<<< \H^1(R, \G) 
                                        @>>> \H^1(K, \G_\eta)
\end{CD}
\]
where the vertical maps are those corresponding to evaluation at $\Ps
\in \Xs(k)$, $P \in \X(R)$ and $P_\eta \in \X_\eta(K)$ respectively.
\end{lemma}
\begin{proof}
The left-hand square commutes because, as observed above, pulling back
by the morphism $\Ps\colon  \Spec k \to \Xs$ gives functorial maps on
cohomology.

We now treat the right-hand two squares.  Observe that, if $S$ is
any scheme over $k$, then every $R$-morphism from $S$ to $\G$ factors
through the special fibre $\Gs$, by the universal property of fibre
products.  So the sheaves on $\Xs$ defined by $\G$ and $\Gs$ are the
same.  Similarly, the sheaves on $\X_\eta$ defined by $\G$ and $\G_\eta$ are
the same.  We may therefore replace $\Gs$ or $\G_\eta$ in the diagram by
$\G$ wherever they occur.  The maps in the diagram are then all seen
to be natural morphisms coming from morphisms of schemes, given by the
following diagram.
\[
\begin{CD}
\Xs     @>>> \X @<<<      \X_\eta \\
@AAA         @AAA         @AAA \\
\Spec k @>>> \Spec R @<<< \Spec K
\end{CD}
\]
Since this diagram commutes, so does the one above.
\end{proof}

\begin{proof}[Proof of Theorem~\ref{thm:main}]
We begin by replacing the model $\X$ with its maximal smooth
subscheme, which means simply removing the singular locus of the
special fibre.  Now $X(K)_\sm$ consists of all those points of $X(K)$
which still extend to $\X(R)$.

By Theorem~\ref{thm:tech}, the cohomology class of the torsor $Y\to X$
extends to a cohomology class $\alpha \in \H^1(\X,\G)$.  There is
therefore a diagram
\[
\begin{CD}
\Xs(k)      @<<< \X(R)      @>>> X(K)_\sm \\
@VVV             @V{\alpha}VV    @VV{Y}V \\
\H^1(k,\Gs) @<<< \H^1(R,\G) @>>> \H^1(K,G)
\end{CD}
\]
in which the right-hand vertical map comes from the torsor $Y$; the
middle one from the class $\alpha$; and the left-hand one from the
image of $\alpha$ in $\H^1(\Xs, \Gs)$.  This diagram commutes by
Lemma~\ref{lem:eval}.  By construction, the map $\X(R) \to X(K)_\sm$
is surjective; since $\X$ is separated, it is bijective.  The maps in
the bottom row are both isomorphisms, by Lemmas~\ref{lem:const}
and~\ref{lem:hensel}, hence the result.
\end{proof}

\begin{proof}[Proof of Corollary~\ref{cor:finite}]
If $k$ is finite, let $\beta$ be the class in $\H^1(\Xs, \Phi(\G))$ which
is the image of $\alpha$ under the map $\H^1(\Xs, \Gs) \to \H^1(\Xs,
\Phi(\G))$.  Then we can extend the above diagram:
\[
\begin{CD}
\Xs(k)      @=    \Xs(k)      @<<< \X(R)      @>>> X(K)_\sm \\
@V{\beta}VV        @VVV        @V{\alpha}VV       @VV{Y}V \\
\H^1(k,\Phi(\G)) @<<< \H^1(k,\Gs) @<<< \H^1(R,\G) @>>> \H^1(K,G)
\end{CD} \text{.}
\]
The diagram still commutes, by Lemma~\ref{lem:eval}.  The maps in the
bottom row are all isomorphisms (by Lemma~\ref{lem:finite}), and the
result again follows.
\end{proof}

\section{Examples and applications}\label{sec:examples}

The use of Theorem~\ref{thm:main} is that it allows us to reduce
questions about the arithmetic of $X$ to questions about the geometry
of the special fibre $\X$.  In this section, we demonstrate this with
some results on torsors under tori.  For these examples, we take $K$
to be a $p$-adic field, so that the residue field $k$ is always
finite.

\begin{proposition}\label{prop:rational}
Let $X$ be a smooth projective variety over $K$, and suppose that $X$
has a model $\X$ such that the smooth locus of the special fibre $\Xs$
is geometrically simply connected. Then, for any $K$-torus $T$ split
by a tame extension of $K$, and any $X$-torsor under $T$, the
evaluation map $X(K) \to \H^1(K,T)$ is constant on $X(K)_\sm$.
\end{proposition}
\begin{proof}
Let $U$ be the smooth subscheme of the special fibre $\Xs$.  Let
$\Phi$ be the component group of the N\'eron model of $T$; in view of
Theorem~\ref{thm:main}, it suffices to show that every class in
$\H^1(U,\Phi)$ is constant.  After base change to $\kb$, $\Phi$
becomes isomorphic to a direct sum of groups each of the form $\Z/n\Z$
or $\Z$.  But, since $\bar{U}$ is simply connected, $\H^1(\bar{U}, \Z)
= \H^1(\bar{U}, \Z/n\Z) = 0$, and therefore $\H^1(\bar{U}, \Phi) = 0$.
Now the Hochschild--Serre spectral sequence gives an isomorphism
$\H^1(k, \Phi) \to \H^1(U, \Phi)$, completing the proof.
\end{proof}

In particular, let $X$ be a smooth del Pezzo surface; then
the special fibre of $\X$ may be a del Pezzo surface with isolated
singularities.  Proposition~\ref{prop:rational} applies if the
singularities are of certain types, as listed in~\cite[\S
  3]{KM:MAMS-1999}.  These are precisely the singularity types such
that $\Pic \bar{U}$ is free; this can easily be verified by computing
$\Pic \bar{U}$ as the quotient of the relevant root lattice by the
sublattice generated by the exceptional curves of a resolution.

\begin{remark}
We can also make conclusions about del Pezzo surfaces with other
reduction types.  For example, suppose that $X$ is a smooth del Pezzo
surface, and $\Xs$ is the cone over a smooth curve $C$ of genus $1$.
After removing the vertex to obtain the smooth locus $U$, there is an
isomorphism $\H^1(U,\Phi) \cong \H^1(C,\Phi)$.  Evaluating $X$-torsors
under $T$ becomes a question about $C$-torsors under $\Phi$.  For an
application of this in terms of the Brauer--Manin obstruction,
see~\cite{Bright:cubics}.
\end{remark}

\begin{remark}
Applying the techniques of this article can only give information
about the points of $X$ reducing to the smooth locus of $\Xs$.  For
more complete information, we should first construct a weak N\'eron
model for $X$, which is a model such that any $K$-point of $X$ extends
to a smooth $R$-point of (some component of) the model.
\end{remark}

Following the work of Colliot-Th\'el\`ene and Sansuc~\cite{CTS:DMJ-1987}
on torsors under tori, we can also use information about the Galois
action on the Picard group of $X$ to deduce results about $X$-torsors
under tori.

\begin{theorem} \label{thm:tamegal}
Let $X$ be a smooth, geometrically integral, projective variety over a
$p$-adic field $K$, and suppose that the Galois action on the
geometric Picard group of $X$ is tame.  Let $\X/R$ be a model of $X$.
Let $T$ be \emph{any} torus over $K$, and let $Y \to X$ be an
$X$-torsor under $T$.  Let $\T$ denote the N\'eron model of $T$.  Then
the conclusion of Corollary~\ref{cor:finite} holds, i.e.\ the
evaluation map $X(K)_\sm \to \H^1(K,T)$ coming from $Y$ factors
through the $k$-points of the special fibre $\Xs$, as follows:
\[
\begin{CD}
X(K)_\sm @>>> \Xs(k) \\ @V{Y}VV @VV{g}V \\ \H^1(K,T) @>{\cong}>>
\H^1(k,\Phi(\T))
\end{CD}
\]
where the map $g$ comes from an element of $\H^1(\Xs, \Phi(\T))$.  The
bottom isomorphism is the composite of the isomorphisms
\[
\H^1(K,T) \leftarrow \H^1(R,\T) \to \H^1(k, \T_s) \to \H^1(k, \Phi(\T))
\]
of Section~\ref{sec:constant}.
\end{theorem}
\begin{proof}
As in the proof of
Theorem~\ref{thm:main}, we first turn $\X$ into a smooth model by
removing the singular subscheme of the special fibre.  Now, although
it is not necessarily true that $\R^1_\sm \eta_* T = 0$, we will show
that $Y$ nevertheless extends to a cohomology class in $\H^1(\X,\T)$;
the remainder of the proof is exactly as for Theorem~\ref{thm:main}.

Let us first define a sub-torus $S \subseteq T$, as follows.  Let $M =
\Hom_\Ksep(T, \Gm)$ be the module of characters of $T$; it is a finitely
generated, free Abelian group with an action of $\Gamma =
\Gal(\Ksep/K)$.  Let $\Gamma_1 \subseteq \Gamma$ be the wild inertia
group.  Let $N$ be the largest free quotient of $M$ on which
$\Gamma_1$ acts trivially.  Let $S$ be the sub-torus of $T$ dual to
$N$; then $S$ is the largest sub-torus of $T$ which is split by a
tamely ramified extension of $K$.

Now we will show that $\H^1(X, T)$ is generated by the images of
$\H^1(X, S)$ and $\H^1(K, T)$.  We may assume that $X(K)$ is non-empty,
since otherwise the conclusion of the theorem is vacuous.  Write $\Xb
= X \times \Ksep$; Colliot-Th\'el\`ene and
Sansuc~\cite[Lemme]{CTS:CR-1976} have shown that there is an exact
sequence
\[
0 \to \H^1(K, T) \to \H^1(X, T) \to \Hom_K(M, \Pic\Xb) \to 0
\]
which is functorial in $T$.  There is therefore a commutative diagram
\[
\begin{CD}
0 @>>> \H^1(K,S) @>>> \H^1(X,S) @>>> \Hom_K(N, \Pic\Xb) @>>> 0 \\
@.     @VVV          @VVV          @VV{\cong}V              \\
0 @>>> \H^1(K,T) @>>> \H^1(X,T) @>>> \Hom_K(M, \Pic\Xb) @>>> 0
\end{CD}
\]
where, since $\Gamma_1$ acts trivially on $\Pic \Xb$, the right-hand
vertical arrow is an isomorphism.  It follows that $\H^1(K,T)$ surjects
onto the cokernel of $\H^1(X,S) \to \H^1(X,T)$; therefore $\H^1(X,T)$ is
generated by the images of $\H^1(X,S)$ and $\H^1(K,T)$.

The image of $\H^1(\X, \T) \to \H^1(X,T)$ contains both $\H^1(X,S)$
(since $S$ is split by a tamely ramified extension) and $\H^1(K,T)$
(by Lemma~\ref{lem:const}).  It therefore is the whole of $\H^1(X,T)$.
\end{proof}

As an application, we deduce the following corollary.  The result is
already known: for example, it follows from the result of
Colliot-Th\'el\`ene~\cite{CT:IM-1983} that the Chow group of 0-cycles
on such a surface is trivial.

\begin{corollary}
Let $X$ be a smooth, projective, geometrically rational surface over
$K$ and suppose that there exists a smooth, proper model of $X$ over
$R$.  Let $T$ be any torus over $K$, and let $Y \to X$ be an
$X$-torsor under $T$.  Then the evaluation map $X(K) \to \H^1(K,T)$
defined by $Y$ is constant.
\end{corollary}
\begin{proof}
In this case of good reduction, $\Pic \Xb$ is unramified as a Galois
module, as follows.  Let $\X$ be a smooth, proper model of $X$ over
$R$.  Write $\Xsb = \X \times \kb$, and recall that there is a
reduction map from $\Pic \Xb$ to $\Pic \Xsb$, which under our
hypotheses is an isomorphism~\cite[Proposition~3.4.2]{Harari:DMJ-1994}
and clearly respects the Galois action.  Since the inertia group $I$
acts trivially on $\Pic\Xsb$, we deduce that $I$ also acts trivially
on $\Pic\Xb$.  So Theorem~\ref{thm:tamegal} applies.

The same argument as in the proof of Proposition~\ref{prop:rational}
now shows that the evaluation map defined by $Y$ is constant.
\end{proof}

\paragraph{Acknowledgements}
I thank the anonymous referees for comments leading to improvements in
the presentation of this article.

\bibliographystyle{abbrv}
\bibliography{martin}

\end{document}